# A Macro-Micro Approach to Reconstructing Vehicle Trajectories on Multi-Lane Freeways with Lane Changing

Xuejian Chen [a, b], Guoyang Qin [a], Toru Seo [b], Ye Tian [a], Jian Sun [a*]

[a] *College of Transportation Engineering & Key Laboratory of Road and Traffic Engineering, Ministry of Education, Tongji University, Shanghai 201804, China*
[b] *Department of Civil and Environmental Engineering, Tokyo Institute of Technology, Tokyo 152-8552, Japan*

**Abstract**

Vehicle trajectories can offer the most precise and detailed depiction of traffic flow and serve as a critical component in traffic management and control applications. Various technologies have been applied to reconstruct vehicle trajectories from sparse fixed and mobile detection data. However, existing methods predominantly concentrate on single-lane scenarios and neglect lane-changing (LC) behaviors that occur across multiple lanes, which limit their applicability in practical traffic systems. To address this research gap, we propose a macro-micro approach for reconstructing complete vehicle trajectories on multi-lane freeways, wherein the macro traffic state information and micro driving models are integrated to overcome the restrictions imposed by lane boundary. Particularly, the macroscopic velocity contour maps are established for each lane to regulate the movement of vehicle platoons, meanwhile the velocity difference between adjacent lanes provide valuable criteria for guiding LC behaviors. Simultaneously, the car-following models are extended from micro perspective to supply lane-based candidate trajectories and define the plausible range for LC positions. Later, a two-stage trajectory fusion algorithm is proposed to jointly infer both the car-following and LC behaviors, in which the optimal LC positions is identified and candidate trajectories are adjusted according to their weights. The proposed framework was evaluated using NGSIM dataset, and the results indicated a remarkable enhancement in both the accuracy and smoothness of reconstructed trajectories, with performance indicators reduced by over 30% compared to two representative reconstruction methods. Furthermore, the reconstruction process effectively reproduced LC behaviors across contiguous lanes, adding to the framework's comprehensiveness and realism.

## 1. Introduction

Vehicle trajectories are a valuable source of traffic information in both spatial and temporal domains. They have been widely used in various traffic-related applications, such as traffic state estimation, traffic flow modeling, signal optimization, and energy emission estimation (Li et al., 2020). The availability of fully sampled vehicle trajectories is vital for a comprehensive description of the entire traffic flow and underpins its diverse applications. However, obtaining fully sampled vehicle trajectories through direct observation in the natural world, such as using video cameras (NGSIM, 2006; Seo et al., 2020) or unmanned aerial vehicles (Barmpounakis and Geroliminis, 2020; Krajewski et al., 2018), can be laborious and costly. Furthermore, the current mainstream traffic sensors, including both fixed and mobile sensors, can only provide partial trajectory data due to their low deployment and low penetration rates.

Therefore, reconstructing vehicle trajectories from limited detected traffic data has become an increasingly important research topic in recent years. Various techniques have been devised to reconstruct trajectories based on prevailing data sources. Wang et al. (2020) estimated the trajectories of all human-driven vehicles within mixed traffic flows under a connected vehicle environment. Tsanakas et al. (2022) generated virtual vehicle trajectories for emission estimation under a multi-source data environment. In our previous work (Chen et al., 2022a), we fused fixed sensor and probe vehicle (PV) trajectories to reconstruct individual trajectories of undetected vehicles for freeways. However, these methods are limited to single-lane traffic scenarios, leaving the effect of lane-changing (LC) behaviors to be explored.

*Corresponding author, email: sunjian@tongji.edu.cn

Existing research has achieved desirable accuracy in trajectory reconstruction on single lanes. However, extending the single-lane methodologies to multi-lane scenarios is a non-trivial task because the two scenarios differ in several essential aspects. There are three main challenges associated with this extension: *(1) Trajectory degree of freedom is relaxed from solely CF to jointly CF and LC behaviors.* In the single-lane scenario as illustrated in Fig. 1(a), driving interactions between consecutive vehicles are primarily governed by the car-following dynamics, and the lane boundaries suppress any intention for LC. Consequently, the trajectory adjustment of the ego vehicle in blue is simply confined by its preceding and succeeding vehicles. In contrast, the multi-lane scenarios allow for more versatility in the trajectory of the ego vehicle since the suppression due to lane boundaries is lifted and the vehicle is enabled to decide whether to follow its leader or change lanes; *(2) Traffic state constraints are expanded for LC-enabled trajectories.* Macroscopic traffic state emerges from a collection of microscopic trajectories and governs the geometry of each individual microscopic trajectory. In the multi-lane scenario, the governing law becomes complex as the velocity difference between adjacent lanes becomes a new force to guide the trajectories in addition to the relatively simple longitudinal force in single-lane scenarios. To this end, the motions of LC vehicles need to obey the traffic state constraints on both lanes; *(3) Surrounding vehicle trajectories are sensitive to the LC positions inferred.* Inferring LC positions in a sparse sensing environment is a hard problem due to significant information loss. This can be seen in Fig. 1(b), where given the observed boundary trajectories in red, two possible reconstruction plans (trajectories in dashed lines) can be created for the two candidate LC positions of the ego vehicle in blue lines. Determining which construction plan is more likely to be close to the ground truth is the core challenge. It is noted that the candidate position Y causes more sensitive fluctuations of surrounding vehicles' trajectories than the position X. This observation can be leveraged to reduce the uncertainty of LC position inference.

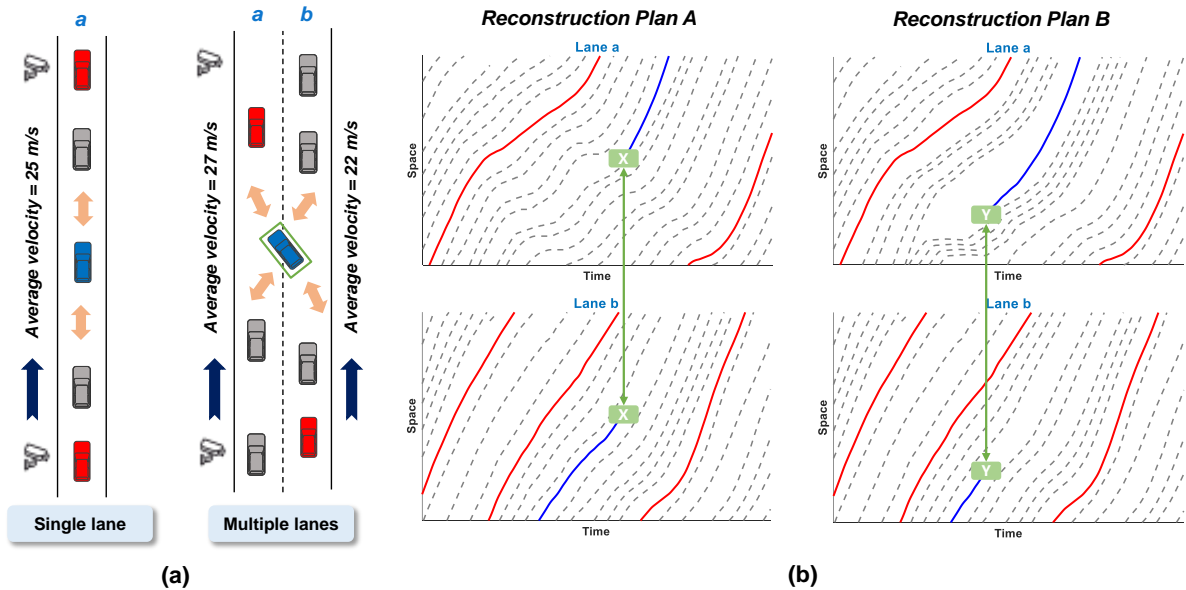

**Fig. 1. Illustration of multi-lane trajectory reconstruction challenges: (a) Comparison between**

**single-lane scenario and multi-lane scenario: (b) Impact of LC position for trajectory accuracy.**

To address the associated challenges, we propose a multi-lane trajectory reconstruction approach to accommodate lane-change behavior under sparse fixed and mobile sensing. Our major contributions are outlined below:

1) We propose a macro-micro framework to decouple the multi-lane task into two components: LC position identification and lane-based individual trajectory reconstruction. This approach overcomes the lane boundary and improves computational efficiency.
2) To identify the plausible LC positions, we derive velocity contour maps from fixed and mobile sensor data. The comprehensive traffic state provides constraints for the motions of collective vehicles, while velocity difference across multiple lanes emerges as a significant factor to guide LC behaviors.
3) To reconstruct lane-based individual trajectories, we develop four trajectory estimation algorithms based on extended CF models. These candidate trajectories offer precise baselines for each undetected vehicle, and determine the potential range of LC positions that satisfies safety distances.
4) To jointly infer CF and LC behaviors, we propose a two-stage trajectory fusion algorithm to fuse candidate trajectories under the constraints of the established velocity contour map, meanwhile optimizing LC position by balancing macro traffic state variations and micro trajectory fluctuation.
5) We evaluate the proposed method using field-surveyed datasets. The results demonstrate significant improvements in reconstruction accuracy and smoothness, as well as more effectively capturing the LC behaviors compared to two non-integrated methods.

The structure of the remaining paper is organized as follows. We first review the related literature about vehicle trajectory reconstruction in Section 2, then describe the problem of trajectory reconstruction and the detection environment in Section 3. Then we propose the framework and outline the three main reconstruction algorithms in Section 4. Finally, we conduct the evaluation using NGSIM dataset in Section 5 and draw the conclusion and future research directions in Section 6.

## 2. Literature Review

Vehicular trajectory reconstruction has garnered significant attention over the last few years. Several techniques have been devised to reconstruct vehicular trajectories using diverse data sources. In this section, existing trajectory reconstruction studies are classified based on their data sources. This classification can facilitate useful comparisons among different techniques for evaluating their effectiveness.

The first category pertains to fixed sensors, which have already been employed for traffic state estimation (Seo et al., 2017) and path flow reconstruction (Rao et al., 2018). However, the utilization of this single detection data alone for achieving precise vehicle trajectory reconstruction has received relatively little scholarly scrutiny. Coifman (2002) utilized a simplified kinematic wave theory to reconstruct trajectories for freeways using loop detectors but overlooked the variability of vehicle velocities. Drawing inspiration from Coifman's work, subsequent scholars refined the approach by estimating individual vehicle velocities from upstream and downstream fixed sensors, and then reconstructing trajectories based on kinematic equations (Chen et al., 2014; Van Lint, 2010) However, these methods tend to ignore sophisticated vehicle behaviors like CF and acceleration/deceleration processes. Although they have demonstrated favorable performance in estimating travel times, the resulting trajectories were unsuitable for more granular applications, such as energy emissions and traffic flow oscillation analysis.

The rapid advancements in PV and connected vehicle technologies have brought mobile sensor-based methods into the forefront of trajectory reconstruction research. While these methods overcome the limitations of fixed sensors in terms of detection range, they are hindered by the low upload frequency and low penetration rates. To address the above issues, several traffic flow models, such as Kinematic Wave model or CF, have been utilized to reconstruct fully-sampled vehicle trajectories. However, Kinematic Wave model-based methods exhibit limitations in incorporating the stochastic volatility of velocities and capturing realistic acceleration patterns, resulting in reconstructed trajectories that conform strictly to a constant acceleration profile (Chen et al., 2020; Sun and Ban, 2013). While CF methods can approximate acceleration or deceleration, they necessitated a higher penetration rate due to the accumulation of estimated errors (Goodall et al., 2016). Furthermore, these mobile sensor-based methods may encounter difficulties when applied in multi-lane environments due to the neglect of LC behaviors.

Nowadays, certain scholars are engaged in reconstructing trajectories in the context of connected and automated vehicles. Compared to traditional PVs, these vehicles are capable of collecting not only their own trajectories but also those of surrounding vehicles within a certain range (Seo and Kusakabe, 2015). Wang et al. (2020) and Chen et al. (2021) have respectively employed the Wiedemann Model and the Intelligent Driver Model to estimate the trajectories within mixed traffic flows. Meanwhile, Qi and Chen (2021) applied Bayesian network methods to infer the trajectories of multiple vehicles located ahead. Chen et al. (2022b) have combined shockwave and CF models to generate fully-sampled trajectories for signalized intersections. However, it should be noted that such vehicles remain at an extremely low penetration rate currently, and they can only collect specific trajectories within a limited period. Therefore, their approaches may be considered inadequate when implemented in PV data environments because of the high data collection requirements, and they can only be feasibly employed in single-lane situations.

In light of the aforementioned limitations of both fixed and mobile sensors, multi-source data fusion has emerged as a viable alternative and garnered increased attention in recent years. The fusion of heterogeneous data has yielded improvements in various applications, including traffic state estimation (Deng et al., 2013) and emission consumption (Jiang et al., 2017). Several studies have endeavored to fuse fixed sensors with vehicle IDs (e.g., automatic vehicle identification and video camera) as well as PVs. For instance, Mehran et al. (2012) applied variational theory to reconstruct trajectories at signalized intersections. But the reconstructed trajectories were insufficiently precise, as the primary objective was to estimate travel time. Alternatively, other studies have focused on fusing data from fixed sensors without vehicle IDs, such as loop detectors and radar, with mobile sensors. They have generated macro speed surface from multi-source data, then employed it to estimate virtual vehicle trajectories for freeways (Tsanakas et al., 2022; Van Lint and Hoogendoorn, 2010). Although these methods have achieved good performance in individual trajectory estimations, they were still restricted by the single lane, primarily due to the lack of consideration in the interactive vehicular behaviors between adjacent lanes.

In recent years, there has been a growing recognition of the significance for investing multi-lane scenarios owing to their closer resemblance to real-world traffic systems. In comparison to single-lane scenarios, multi-lane roadways exhibit inherent complexity and instability due to the stochasticity introduced by LC behaviors. However, prevailing research for multi-lane scenarios has predominantly focused on the macro perspective, i.e., estimating various traffic parameters such as flow and density from multi-source data, which aims to serve as a foundation for lane-level traffic management (Bekiaris-Liberis et al., 2017; Kyriacou et al., 2021; Liu et al., 2023). From the micro perspective, Rey et al. (2019) estimated vehicle trajectories with overtaking behaviors by relaxing the first-in-first-out assumption. Nonetheless, the reconstruction encompassed a superposition of trajectories for different lanes, rather than capturing individual lane-based trajectories, which fell short in addressing the interactive behavior among vehicles within distinct lanes. Therefore, reconstructing precise lane-based trajectories still poses a significant challenge in the spatiotemporal sparse data environments, and it is imperative to develop a novel reconstruction method with LC behaviors for multi-lane scenarios.

## 3. Problem Description

The primary aim of the proposed method is to reconstruct vehicle trajectories on multi-lane freeways, and the LC behavior is considered during the reconstruction process, i.e., if a vehicle moves from Lane a to Lane b, its trajectory will be simultaneously removed from Lane a and appear in Lane b, as outlined by green circle in Fig. 2. Two kinds of detected data are involved as input in the studied data environment: fixed sensors located upstream and downstream, and PVs with low penetration rates. Fixed sensors are required to provide individual vehicle information including ID, arriving time and instantaneous velocity, as show in Fig. 2 with blue dots. Such information is accessible under the present technology level, like using automatic vehicle identification systems or roadside unit sensors. Meanwhile PVs supply complete trajectories with lane information, which can be achieved through technologies such as vehicle-to-vehicle or vehicle-to-infrastructure communication (Dey et al., 2016). But they only account for a small proportion of traffic flow, as shown by the red lines. Regarding the typical distance between two consecutive fixed sensors on urban freeway is 500~800 meters, the range of reconstruction is also around 500 meters.

All the grey dash lines in Fig. 2 denote those non-probe vehicles trajectories, which are the reconstruction objects of this study. Some of these vehicles maintain their lane throughout the journey and are detected twice at the same lane, while others change their lanes and are detected in different lanes at upstream and downstream. However, observing LC positions from fixed and mobile sensors is difficult due to the limited observations. Additionally, PVs' trajectories are far less than non-probe vehicles, which also increases the complexity of full-sampled trajectory reconstruction task.

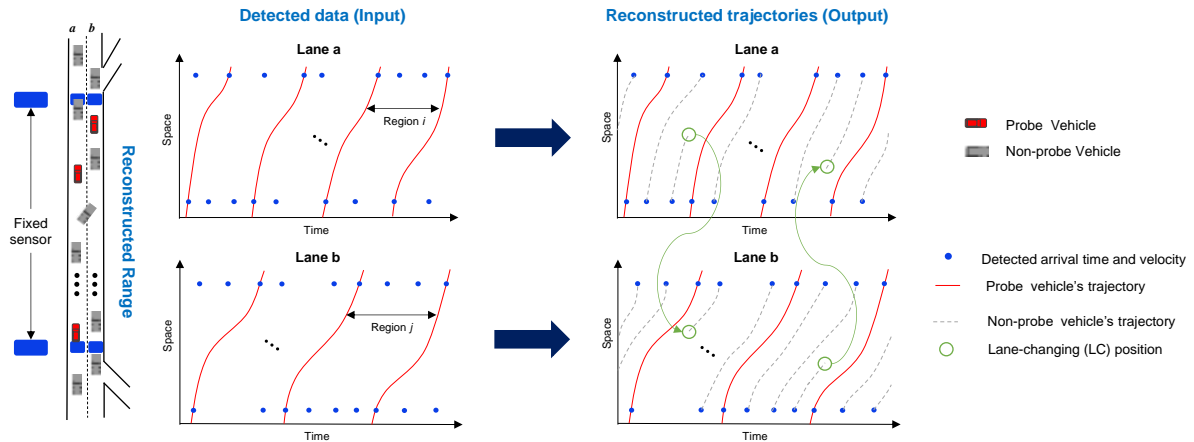

**Fig. 2. The data environment and reconstruction objects in the studied area.**

Instead of reconstructing macroscopic traffic state like density, volume or queue length, the proposed method takes priority to generate second-level trajectories for individual vehicles. The goal of multi-lane trajectory reconstruction can be described as:

$$\min \sum_{l=1}^{L} \sum_{n=1}^{N} \sum_{t=1}^{T} \left| Traj_{re}^{l,n,t} - Traj_{gt}^{l,n,t} \right|$$

s.t.

$$Traj_{re}^{l,n,t} \propto TS_{macro}^{l,t}$$
$$Traj_{re}^{l,n,t} \propto CF_{micro}\left(Traj_{re}^{l,n-1,t}, Traj_{re}^{l,n+1,t}\right)$$
$$Traj_{re}^{1,n',t'} = Traj_{re}^{2,n',t'}, n' \in LC\ vehicle, t' = LC\ time$$

(1)

where $Traj_{re}^{l,n,t}$ and $Traj_{gt}^{l,n,t}$ denote the reconstructed and ground-truth trajectory of undetected vehicle $n$ at time $t$ on lane $l$. Therefore, the objective function is formulated to minimize the difference between reconstructed and ground-truth trajectories. $TS_{macro}^{l,t}$ denotes the traffic state information of lane $l$ at time $t$, and the first constraint enforces that the reconstructed trajectories adhere to the propagation of macroscopic traffic flow. $CF_{micro}$ represents governing principles of car-following behaviors, and the second constraint ensures that the reconstructed trajectories conform to the microscopic motion involving leading and following vehicles. The final constraint guarantees that the lane change position remains consistent between adjacent lanes for those vehicles with LC behaviors. For a comprehensive understanding of the intricacies associated with these constraints, we refer readers to the methodology section.

To simplify the complexity of reconstruction, the space-time diagram is partitioned into multiple regions based on the trajectories of PVs, as illustrated in Fig. 2. Each region is defined as the interval between two consecutive PVs on the same lane, encompassing two PVs and numerous non-probe vehicles. Consequently, PV assumes the absence of lane-changing behavior, which is acceptable given the relatively low penetration rate observed in this study. It should be noted that our focus in the following sections is solely on estimating undetected trajectories within one specific region, with the other regions can be reconstructed in a similar way.

## 4. Methodology

This section first describes the flowchart of the proposed reconstruction approach, which involves three modules: macro-level module, micro-level module and integration module. The overall flowchart is illustrated in Fig. 3. As seen from it, the multi-lane reconstruction task encompasses two primary components, namely LC position identification and lane-based individual trajectory reconstruction. The process of identifying a plausible LC position entails the generation of velocity contour maps within the macro-level module. These maps serve a dual purpose: imposing constraints on the collective movements of vehicles and offering guidance for LC estimation. On the other hand, individual trajectory reconstruction involves the utilization of trajectory estimation algorithms within the micro-level module. These algorithms generate candidate trajectories and determine the potential range of LC positions. Subsequently, a two-stage trajectory fusion algorithm is employed to jointly integrate the macro and micro modules, resulting in the reconstruction of LC-enabled trajectories. The main work of three modules are summarized as follows.

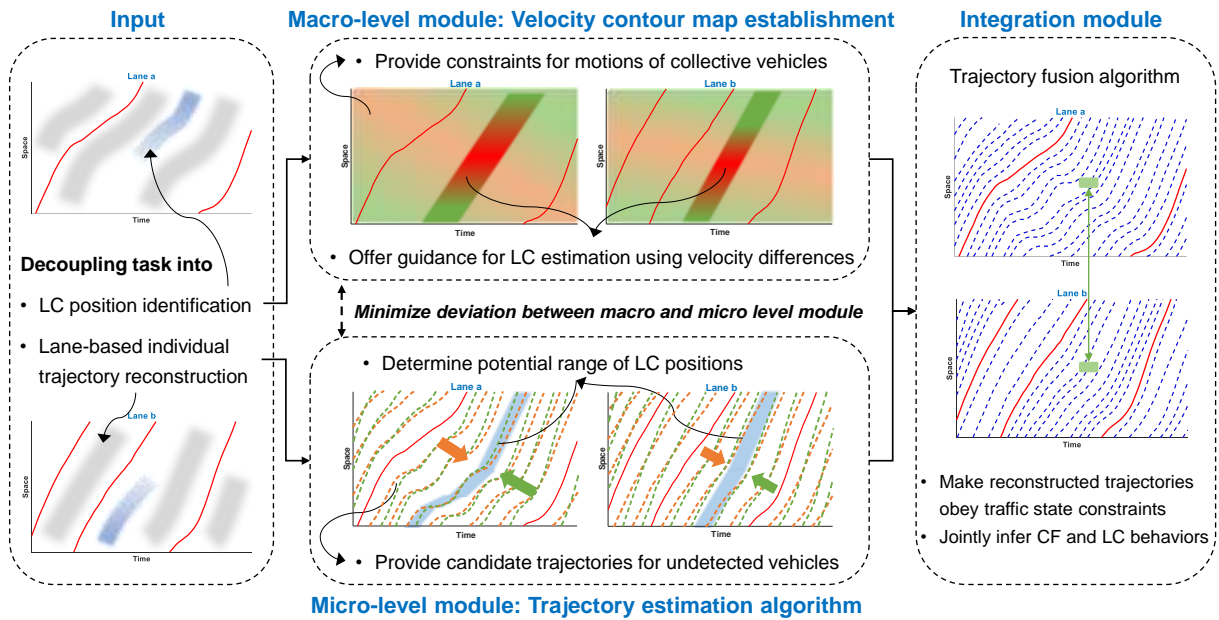

**Fig. 3. Overview of the proposed method for vehicle trajectory reconstruction.**

**Macro-level Module**: Traffic state refers to the collective behavior of vehicles and is typically represented by their trajectories. Therefore, a complete velocity contour map can be used to constrain aggregated vehicle movements and provide crucial information for trajectory reconstruction. Furthermore, velocity difference between adjacent lanes also has a substantial impact on driving behaviors and holds potential in determining LC positions. To obtain comprehensive traffic states, the fixed data and PV trajectories are incorporated to estimate the velocity contour map for entire space-time diagram through a modified adaptive smoothing method (see Section 4.1 for detail).

**Micro-level Module**: Since the velocity contour map can only supply macro traffic parameters, the trajectories reconstructed based on macro-level module ignore the precise CF behaviors and might lead to overlapping problem. To mitigate the issue of the micro overlapping and error accumulation, CF model and its extended Inverse Car-Following (ICF) model are utilized to generate lane-based candidate trajectories from upstream and downstream observations, aiming to cover the possible location ranges of each non-probe vehicle (see Section 4.2 for detail).

**Integration Module**: As the macro traffic state can constrain vehicle movement but fails to capture CF behavior, while micro CF models can provide trajectories but are limited by the neglection of LC behaviors, we proposed a two-stage trajectory fusion algorithm that infers the constraints from both macro traffic state and micro CF and LC behaviors. The first stage is to calculate the weights of candidate trajectories by minimizing the deviation with the velocity contour map; and the second stage is to identify the optimal LC positions based on the traffic state variations (see Section 4.3 for detail).

## 4.1 Space-time Velocity Contour Map Establishment at the Macro Level

The relationship between traffic state and vehicle trajectories is closely intertwined since the aggregate traffic state influences the movement of individual vehicles, and in turn, the collective actions of vehicles also shape the traffic state. Thus, obtaining a comprehensive traffic state is crucial for accurately reconstructing fully-sampled vehicle trajectories. Besides, the traffic state has a significant impact on the driving behaviors, e.g., drivers are more likely to change their lanes in higher-density area in order to seek faster speeds or to avoid congestion (Keyvan-Ekbatani et al., 2016). Hence the LC positions can be plausibly estimated according to the traffic state variation.

However, the complete traffic state cannot be attained directly from both fixed-detected data and PV trajectories due to their low deployment and low penetration rates. Several model-based methods have been proposed to estimate traffic state using partial observations and they can be classified into two categories in general. The first category is based on traffic flow models, like using Lighthill–Whitham–Richards or Aw–Rascle–Zhang models (Seo and Bayen, 2017), and extended filters are usually combined to improve their accuracy. One can refer to a survey by Seo et al. (2017) for more information. But these models require discretizing the space-time diagram into multiple cells, and the data is transformed to aggregate format for further calculation. In our scenario, we collect the information from individual vehicles and demand disaggregate velocity points as output. Hence, these traffic flow-based models are not suitable as they might result in significant loss of accuracy during the transformation from aggregate to disaggregate data.

The second category is based on interpolation models, including linear interpolation and non-parametric interpolation methods. However, as discussed by Tsanakas et al. (2022), linear interpolation models are applicable only in fixed sensor and cannot incorporate data from mobile data source. And sharp changes of speeds and discontinuous acceleration often occur during the interpolation process. On the other hand, non-parametric interpolation models, represented by Adaptive Smoothing Method (ASM), can effectively handle the above issues. ASM, firstly proposed by Treiber and Helbing (2002), has been widely modified and applied to estimate traffic state from heterogeneous data source (Jiang et al., 2017; Van Lint and Hoogendoorn, 2010). The chief novelty of ASM is dividing the propagation of free and congested flows, which makes the estimated state can successfully capture the traffic disturbances and reproduce the stop-and-go waves for urban freeway.

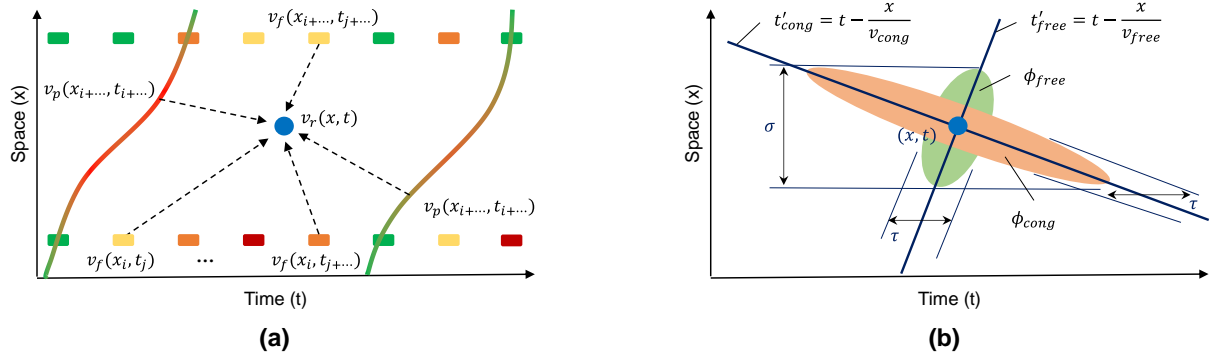

**Fig. 4. Illustrations of macro-level model: (a) disaggregate data detected from fixed sensors and PVs, (b) the anisotropic smoothing kernels for free and congested flows.**

ASM also has a great advantage in dealing with disaggregate data. As seen from Fig. 4(a), the input in our data environment contains the instantaneous velocities detected by fixed sensors (denoted by $v_f(x_i, t_j)$, where $x_i$ and $t_j$ represent the space and time index), and from PVs trajectories (denoted by $v_p(x_i, t_j)$). The object of the macro-level model is to establish the entire velocity contour map, while ASM is particularly well-suited for this task as it can calculate the disaggregate velocity for any point in the space-time diagram. For an arbitrary point $(x, t)$, its velocity $v_r(x, t)$ can be expressed in a discrete convolution form:

$$v_r(x,t) = \frac{\sum_{i=1}^{N}\sum_{j=1}^{M}\phi_0(x-x_i, t-t_j)\cdot V_d(x_i,t_j)}{\sum_{i=1}^{N}\sum_{j=1}^{M}\phi_0(x-x_i, t-t_j)} \quad (2)$$

where $V_d(x_i, t_j)$ denotes the detected velocity at location $x_i$ and time $t_j$. And $\phi_0$ is a kernel function that assigns a weight to each velocity detection $V_d(x_i, t_j)$, with

$$\phi_0(x,t) = \exp(-(\frac{|x|}{\sigma} + \frac{|t|}{\tau})) \quad (3)$$

where the positive constants $\sigma$ and $\tau$ define characteristic widths of the spatial and temporal smoothing. Then two auxiliary velocity surfaces (one for free flow and one for congested flow) are estimated to capture the dynamic feature of traffic flow, and the principal axes are skewed by the anisotropic smoothing kernel using Eq. (4) ~ Eq. (5), as shown in Fig. 4(b).

$$v_r^{free}(x,t) = \frac{\sum_{i=1}^{N}\sum_{j=1}^{M}\beta_{i,j}^{free}(x,t)\cdot V_d(x_i,t_j)}{\sum_{i=1}^{N}\sum_{j=1}^{M}\beta_{i,j}^{free}(x,t)},\ v_r^{cong}(x,t) = \frac{\sum_{i=1}^{N}\sum_{j=1}^{M}\beta_{i,j}^{cong}(x,t)\cdot V_d(x_i,t_j)}{\sum_{i=1}^{N}\sum_{j=1}^{M}\beta_{i,j}^{cong}(x,t)} \quad (4)$$

with

$$\beta_{i,j}^{free}(x,t) = \phi_0(x-x_i, t-t_j - \frac{x-x_i}{v_{free}}),\ \beta_{i,j}^{cong}(x,t) = \phi_0(x-x_i, t-t_j - \frac{x-x_i}{v_{cong}}) \quad (5)$$

where $v_{free}$ and $v_{cong}$ denote the propagation velocities of free flow and congested flow, respectively. The utilization of such anisotropic kernels enables the estimation to account for the characteristic velocity of traffic disruptions. The final velocity contour map is determined based on these two auxiliary velocity surfaces as:

$$v_r^{asm}(x,t) = \omega(x,t)\cdot v_r^{cong}(x,t) + (1-\omega(x,t))\cdot v_r^{free}(x,t) \quad (6)$$

where $\omega(x,t)$ is an adaptive weighting *s*-shape function that depends on the current level of congestion at point $(x, t)$:

$$\omega(x,t) = \frac{1}{2}[1 + tanh(\frac{\hat{v} - min(v_r^{free}(x,t), v_r^{cong}(x,t))}{\Delta v})] \tag{7}$$

where $\hat{v}$ denotes the speed threshold between free and congested flow, while $\Delta v$ denotes the transition width around $\hat{v}$. Readers can refer to Treiber et al. (2011) for more details about the setting of these parameters.

As analyzed above, two types of velocity detections $V_d(x_i, t_j)$ are involved in this study: $v_f(x_i, t_j)$ detected from fixed sensors and $v_p(x_i, t_j)$ collected from PVs. However, the reliability between those data sources varies, e.g., $v_p(x_i, t_j)$ is more accurate in the spatial dimension since PVs' trajectories cover the entire road section, while $v_f(x_i, t_j)$ can efficiently reflect the changes of traffic state at specific locations. Therefore, the source-specific weights should be set to enhance the estimation accuracy in a multi-source data environment. And ASM is modified using Eq. (8).

$$v_r^{asm}(x,t) = \frac{\sum_s \alpha_s \sum_{i,j}^{N,M}[\omega_s(x,t) \cdot \beta_{i,j,s}^{cong}(x,t) + (1 - \omega_s(x,t)) \cdot \beta_{i,j,s}^{free}(x,t)] \cdot v_s(x_i, t_j)}{\sum_s \alpha_s \sum_{i,j}^{N,M}[\omega_s(x,t) \cdot \beta_{i,j,s}^{cong}(x,t) + (1 - \omega_s(x,t)) \cdot \beta_{i,j,s}^{free}(x,t)]} \tag{8}$$

where $s$ denotes the data source and $s \in (v_f, v_p)$, and $\alpha_s$ denotes the corresponding weight of source $s$.

The modified ASM considers the accuracy of different data source as well as the effect of both free and congested flow propagation. Benefitted from that, disaggregate velocities at all spatiotemporal points can be estimated and form a velocity contour map for the entire space-time diagram. This allows for the reconstruction of a comprehensive traffic state, which can then be used to refine the estimated trajectories and determine the plausible LC positions. The details are further explained in the subsequent sections.

### 4.2 Trajectory Estimation Algorithms at the Micro Level

The complete traffic state can be established at macro level, but it fails to generate precise vehicular trajectories since the micro CF behaviors are ignored in the velocity contour map. On the other hand, many CF models are proposed to estimate individual vehicle trajectories and can provide detailed insights into the dynamics of traffic flow. However, conventional CF models are not directly applicable considering the following challenges: *1) Error accumulation.* The low penetration rates lead to a dozen of undetected vehicles that need one-by-one reconstruction based on the leading PV, thus the error would be accumulated to a significant level from downstream to upstream. Besides, the reconstructed trajectories bring potential overlapping with the following vehicles due to the random distributions of PVs through the road; *2) Parameter calibrations.* CF models usually require a large amount of car-following pairs to calibrate their parameters, which limits their usage in our extremely sparse data environment.

The main reason of error accumulation in CF models is that they neglect the impact of the following PV. To overcome the above issues, we extend the conventional CF models in both space and time domains. Specifically, there exists a particular CF interaction between two consecutive vehicles in high-density traffic flow, wherein the speed and distance of the following vehicle vary in response to the movements of its leader. Therefore, the movement of following vehicles can also be utilized to approximately estimate its leading vehicle's trajectory, which is referred to as the Inverse Car-Following model (ICF).

Newell's car-following model (Newell, 1993) is chosen as the basic CF model to mitigate the parameter calibration computation due to its parameters are the same as ASM. Hence a single set of parameters can be used to adapt both the macro and micro models, eliminating the need for separate calibration of the CF models. Moreover, the simple and explicit formulas of Newell model enable easy extension in the space and time domains compared with other CF models such as Intelligent Driver Model and Wiedemann Model. The trajectory of following vehicle $n$ is assumed to be consistent with its leading vehicle $n-1$ in a homogenous space, with a time lag $\eta_n$ and space lag $\theta_n$ shown in Eq. (9).

$$\Delta x_n(t) = \eta_n \cdot v_n(t) + \theta_n \quad (9)$$

where $\Delta x_n(t)$ is the space headway and $v_n(t)$ is the velocity of vehicle $n$ at time $t$. The time lag $\eta_n$ and space lag $\theta_n$ are constant for a given vehicle but may vary across different vehicles, which are timely adjusted based on the detected time headway from fixed sensors. Notably, these two parameters can be written using the velocity of congested flow $v_{cong}$ and jam density $k_j$, as shown in Eq. (10).

$$\eta_n = \frac{1}{v_{cong} \cdot k_j} \text{ and } \theta_n = \frac{1}{k_j} \quad (10)$$

As illustrated in Fig. 5, each non-probe vehicle owns two detected points at downstream and upstream respectively, and the spatiotemporal correlations between detected points and PVs can be classified into four categories: the following vehicle detected by upstream fixed sensor; the following vehicle detected by downstream fixed sensor; the leading vehicle detected by upstream fixed sensor; the leading vehicle detected by downstream fixed sensor. Therefore, four corresponding trajectory estimation algorithms are proposed to handle the above categories and generate candidate trajectories for those non-probe vehicles. The details are as follows.

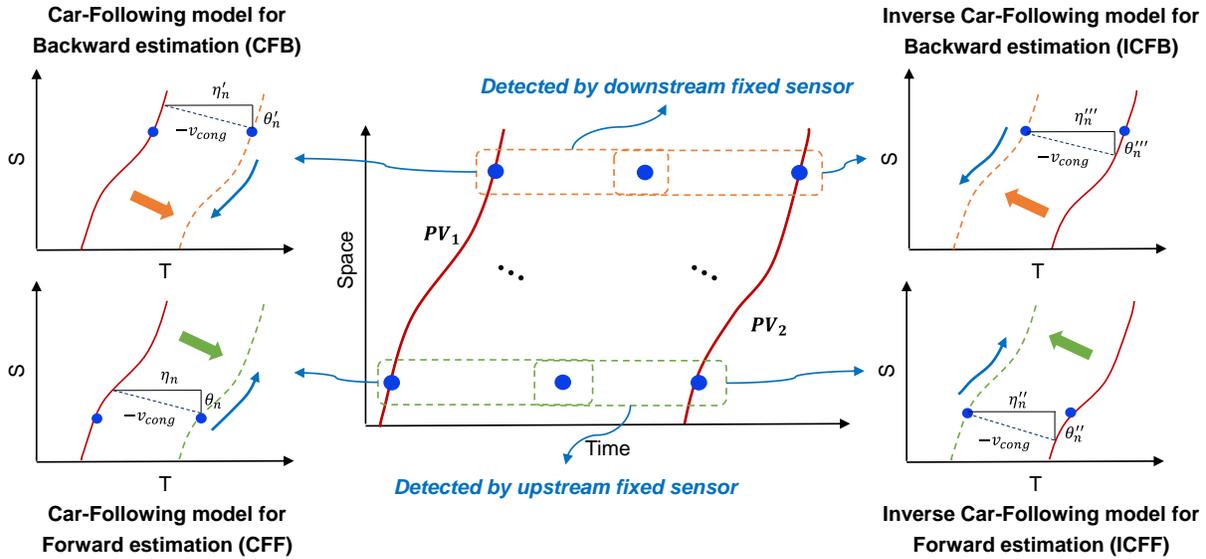

Fig. 5. Four trajectory estimation algorithms proposed to generate candidate trajectories.

For the following vehicle $n$ detected by upstream fixed sensor, CF model is directly applied to generate the undetected trajectories from timestamp $t$ to timestamp $t + 1$ as below, which is named Car-Following model for Forward estimation (CFF).

$$Traj_n^{t+1} = Traj_{n-1}^{t+1-\eta_n} - \theta_n \quad (11)$$

On the other hand, for the following vehicle n detected by downstream fixed sensor with time lag $\eta_n'$ and space lag $\theta_n'$, Car-Following model for Backward estimation (CFB) should be used to generate the trajectories from timestamp $t$ to timestamp $t - 1$:

$$Traj_n^{t-1} = Traj_{n-1}^{t-1-\eta_n'} - \theta_n' \quad (12)$$

Similarly, CF model is easily extended to estimate the leading trajectories by shifting the following trajectory of vehicle $n + 1$. And the Inverse Car-Following model for Forward (ICFF) and Backward estimations (ICFB) can be expressed as follows:

$$\begin{aligned} Traj_n^{t+1} &= Traj_{n+1}^{t+1+\eta_n''} + \theta_n'' \\ Traj_n^{t-1} &= Traj_{n+1}^{t-1+\eta_n'''} + \theta_n''' \end{aligned} \quad (13)$$

where $\eta_n''$ and $\theta_n''$ are the time lag and space lag of the leading vehicle detected at upstream, while $\eta_n'''$ and $\theta_n'''$ are time and space lags detected at downstream. These parameters can be easily calculated using $v_{cong}$ and detected time interval of fixed sensors, as shown in Fig. 5.

With the above analysis, four trajectory estimation algorithms are established based on CF and its extend ICF models, which fully utilize the spatiotemporal correlations of detected information between fixed sensors and PVs. Benefitted from that, each non-probe vehicle owns four candidate trajectories and they can approximately cover all potential locations. However, the precision of these trajectories varies and depends on the estimated vehicle order, i.e., the ground-truth trajectories show more similarities with $PV_1$ at the downstream, suggesting that the CF model is more accurate for front vehicles. Conversely, the ICF model becomes more reliable for the rear vehicles as those trajectories are closer to $PV_2$. Therefore, how to determine the precision degree of the trajectory estimation algorithms and fuse four candidate trajectories plays an important role in reconstruction precision.

### 4.3 Trajectory Fusion Algorithm for Integrated Module

Comprehensive traffic state is facilitated by the modified ASM, which can supply essential information for the movement of vehicle platoons. But the macro level model fails to capture precise CF behaviors and is limited by lane restrictions. On the other hand, trajectory estimation algorithms can provide four candidate trajectories for each non-probe vehicle based on CF and ICF models. However, those trajectories vary in accuracy and their weights need to be determined for further fusion. Moreover, the estimation of LC behavior is also crucial in the realm of traffic management as it can help to optimize traffic flow and reduce the risk of accidents. Nevertheless, LC behaviors is ignored during both the generation of candidate trajectories and the estimation of traffic state, thereby necessitating its inclusion within the following integrated module.

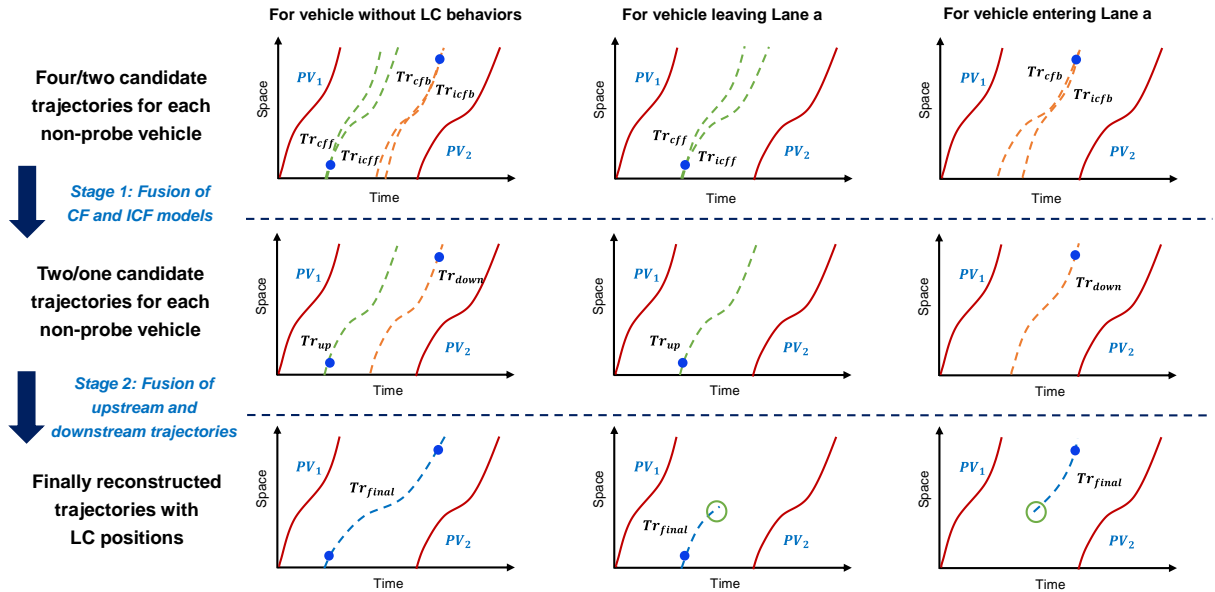

Fig. 6. Two-stage trajectory fusion algorithm for individual vehicles in Lane a.

To solve the above issues, a two-stage trajectory fusion algorithm is proposed to integrate the macro and micro models and jointly infer CF and LC behaviors. The accuracy of four candidate trajectories is evaluated based on the constraint of velocity contour map firstly, and their optimal weights are resolved using dynamic programming. The LC positions are then determined through an optimization-based algorithm that accounts for the variations in traffic state and the safety distance between neighboring vehicles. Therefore, the reconstructed trajectories can effectively balance both micro CF and LC behaviors and macro traffic states. The illustration of the trajectory fusion algorithm is shown in Fig. 6.

As seen from Fig. 6, non-probe vehicles located in Lane a can be classified into three categories in general: *1) Vehicles without LC behaviors.* These vehicles are detected by the upstream and downstream sensors, and own four candidate trajectories within Lane a; *2) Vehicles that leaving from Lane a to Lane b.* Such vehicles are detected solely by the upstream sensor, and possess only two candidate trajectories within Lane a, while the remaining two are generated within Lane b; *3) Vehicles that entering from Lane b to Lane a.* Those are detected exclusively by the downstream sensor, and owning only two candidate trajectories. The proposed trajectory fusion algorithm contains two stages. The first stage is to fuse CF and ICF models under the macro traffic state constraints. This results in a reduction in the number of candidate trajectories for each non-probe vehicle from four/two to two/one, respectively, meanwhile those fused trajectories must cross either the upstream or downstream sensor at the specific detected time. The second stage is to fuse upstream and downstream trajectories for those vehicles without LC behaviors, and estimate LC positions for those vehicles changing lanes. The primary objective is to ensure the consistency of the LC process, i.e., the trajectory should be removed from the original lane and appear in the oriented lane at the same time and position. The details are explained in the following subsections.

*4.3.1 Stage 1: Fusion of CF and ICF models*

The fusion algorithm necessitates the consideration of two factors: designing the moving equations for individual vehicles, and calculating the optimal weights for the whole vehicle platoon. Since the micro trajectory estimation algorithms are established based on Newell models, these candidate trajectories are smooth enough and accord with the kinematic equation. Therefore, the reconstructed position can be easily represented by proportionally adding the CF- and ICF-based candidates. For trajectories of vehicle $n$ estimated at upstream sensor by CFF and ICFF, the fusion equation is formed as:

$$Traj_{up}^n = w_{cff}^n \cdot Traj_{cff}^n + w_{icff}^n \cdot Traj_{icff}^n \tag{14}$$

where $Traj_{cff}^n$ and $Traj_{icff}^n$ denote the candidate trajectories estimated by CFF and ICFF for vehicle $n$ ($n = 1,2, \dots, N$), and $w_{cff}^n$ and $w_{icff}^n$ denote their corresponding weights, respectively. $Traj_{up}^n$ denotes the fused trajectory at upstream. $N$ is the number of non-connected vehicles in each region. The kinematic variables can be calculated as:

$$\begin{aligned} X_{up}^{n,t} &= w_{cf}^n \cdot X_{cff}^{n,t} + w_{icf}^n \cdot X_{icff}^{n,t} \\ V_{up}^{k,n} &= (X_{up}^{k+1,n} - X_{up}^{k-1,n})/(2 \cdot t_{int}) \end{aligned} \tag{15}$$

where $X_{up}^{n,t}$ and $V_{up}^{n,t}$ denote the position and velocity of $Traj_{up}^n$ at time $t$, while $X_{cff}^{n,t}$ and $X_{icff}^{n,t}$ denote the position of $Traj_{cff}^n$ and $Traj_{icff}^n$ at $t$, respectively. $t_{int}$ denotes the time interval.

With the above analysis, the trajectory of individual vehicle can be fused using the moving equations after obtaining the weights of CFF and ICFF. However, the weight calculation must take into account not only individual vehicles but also the whole vehicle platoon. Specifically, the fused trajectory should avoid overlapping with either its leading or following vehicles, while also conforming to the propagation of traffic flow. Therefore, the determination of weights can be formulated as an optimization problem, see Eq. (16).

$$\begin{aligned} \min &\sum_{n=1}^{N} \sum_{t=1}^{T} (V_{up}^{n,t} - V_{asm}^{x,t})^2 \\ \text{s.t.} \quad &V_{up}^{n,t} = f(w_{cff}^n, w_{icff}^n) \\ &w_{cff}^n + w_{icff}^n = 1 \\ &w_{cff}^n > w_{cff}^{n+1} (n = 1,2, \dots, N-1) \end{aligned} \tag{16}$$

where $V_{asm}^{x,t}$ denotes the velocity of the same coordinate point with $V_{up}^{n,t}$, which is estimated from the modified ASM. The objective function is formulated to minimize the velocity deviation between the fused trajectory (Micro-level module) and the established velocity contour map (Macro-level module). In other words, the macroscopic velocity is used as a constraint to aid in determining the optimal weights. $f$ denotes the formulas of moving equations of Eq. (14) ~Eq. (15). To ensure that the fused trajectories reach the upstream sensor at the detected time, the sum of $w_{cff}^n$ and $w_{icff}^n$ should equal 1.

The final constraint specifies that the weight of CFF should decrease monotonically with the reconstructed vehicle order. This constraint is in accordance with the propagation of traffic waves as front vehicles are more influenced by CFF, while ICFF becomes more accurate for rear vehicles. Therefore, the weight proportion of the trajectory estimated by CFF should decrease with increasing vehicle order n. Moreover, the monotonically decreasing weights is also essential for preventing overlap of fused trajectories, particularly as the shape of reconstructed trajectories gradually transition from $PV_1$ to $PV_2$ in the same region. This vehicle platoon optimization can be well-solved using dynamic programming as it is a typical directed and acyclic graph problem. One can refer to Sun and Ban (2013) for more details about the directed graph construction and cost design. The fusion process of CFB and ICFB is akin to that of CFF and ICFF, and the specifics of the process are omitted here due to space limitations.

*4.3.2 Stage 2: Fusion of upstream and downstream trajectories*

The main works in this subsection are bifurcated into two parts: for vehicles without LC behaviors, the estimated trajectories from upstream and downstream sensors are fused once more, and the finally reconstructed trajectories will arrive and depart from the road section at their detected times. For vehicles that do exhibit LC behaviors, the LC positions are initially calculated, followed by the splitting of the trajectories, as depicted in Fig. 6.

The fusion equation for upstream and downstream trajectories is formed using Eq. (17). Where $X_{re}^{n,t}$ denotes the position of finally reconstructed trajectories of vehicle $n$ at time $t$, and $T$ denotes the total time step of reconstructed trajectories. The weight of $Traj_{up}^{n,t}$ is greater than $Traj_{down}^{n,t}$ for the first half, and flip the weight for the second half. Benefitted from the time-varying fusion equation, the reconstructed trajectories are sufficiently smooth to traverse their detected points.

$$X_{re}^{n,t} = \left(\frac{t}{T}\right)^2 \cdot X_{down}^{n,t} + \left[1 - \left(\frac{t}{T}\right)^2\right] \cdot X_{up}^{n,t} \tag{17}$$

For vehicles with LC behaviors, the determination of their changing positions is reliant on the analysis of traffic state variations, given the significant impact of traffic state on driving behaviors as evidenced in prior studies (Xie et al., 2022). Besides, numerous LC models have been developed based on the velocity differences between adjacent lanes (Pang et al., 2020), which further substantiates the rationale of the proposed estimation method.

The estimation of LC behavior is illustrated in Fig. 7. Firstly, the common time steps between two candidate trajectories from different lanes are extracted to generate a feasible area, as shown with the starting time $t_s$ and ending time $t_e$ in Fig. 7(a). This implies that plausible LC behaviors can only occur within the feasible area, and all coordinate points in the feasible area are considered as potential LC positions. Secondly, the estimated LC position is regard as a new "detected sensor", and two candidate trajectories should be adjusted to reach the new detected point, as shown in Fig. 7(a) with extent $Dis$. The adjustment method can refer to the vehicles without LC behaviors as they both own two detected points at upstream and downstream, respectively.

With the above analysis, the optimal LC position is then calculated based on an optimization algorithm that balances the impact of traffic state variations and the adjusted extent of candidate trajectories, while the safe distance between reconstructed trajectories and their neighboring vehicles' trajectories should also be guaranteed during the optimization process. The objective function and corresponding constraints are defined using Eq. (18).

$$\max \frac{|V_{asm}^{x,t}(l_a) - V_{asm}^{x,t}(l_b)| + v_{threshold}}{Dis_n^t + Dis_{threshold}}$$
$$\text{s.t.} \quad t \in [t_s, t_e]$$
$$Dis_n^t = 0.5 \cdot (|X_{up}^{n,t} - X_{down}^{n,t}|)$$
$$L_n^{n-1} > L_{safe} \, \& \, L_n^{n+1} > L_{safe}$$
(18)

where $V_{asm}^{x,t}(l_a)$ and $V_{asm}^{x,t}(l_b)$ denote the velocity estimated from modified ASM for Lane a and Lane b, respectively. $Dis_n^t$ denotes the adjusted extent of two candidate trajectories $Traj_{up}^{n,t}$ and $Traj_{down}^{n,t}$, which is calculated by the absolute half difference of $X_{up}^{n,t}$ and $X_{down}^{n,t}$. $L_n^{n-1}$ and $L_n^{n+1}$ denote the headway between its leading and following vehicles, while $L_{safe}$ denotes the minimum safe distance.

The objective function is formulated to maximize the quotient of the velocity difference and the adjusted extent. Notably, the lower adjusted extent is more amenable to being chosen as LC positions to alleviate abrupt changes in the fused trajectories, while greater velocity difference also holds greater priority owing to the association between traffic state variation and LC behaviors. Hence, the optimized quotient adeptly harmonizes these two factors and yields the optimal LC position. The small values assigned to $v_{threshold}$ and $Dis_{threshold}$ serve to prevent the occurrence of zero numerator and denominator. The first constraint is to ensure that the inferred LC position falls within the feasible area, and the last constraint is to prevent trajectory overlap across multiple lanes. This optimization algorithm can be easily solved due to the restricted scope of the feasible area, as shown in Fig. 7(b).

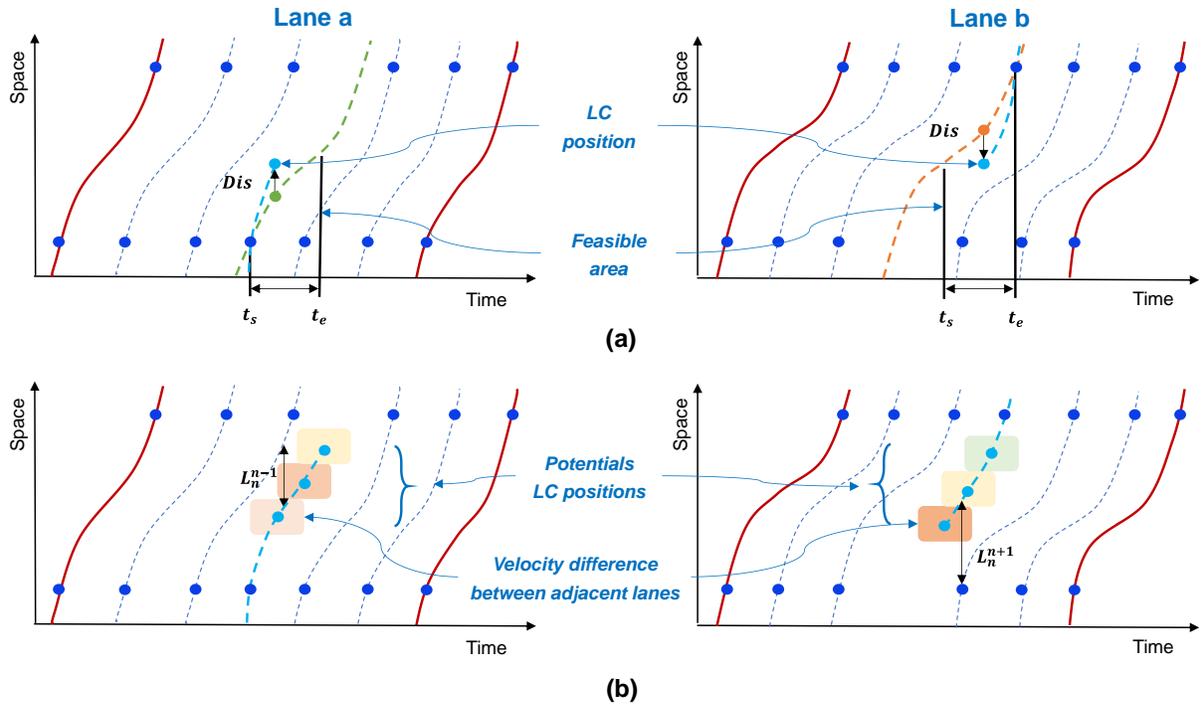

**Fig. 7. Illustration of LC behavior estimation: (a) The adjusted extent of two candidate trajectories; (b) The constraints to determine LC positions.**

In summary, the proposed two-stage trajectory fusion algorithm can integrate the macro and micro models. And the reconstructed trajectories not only conform to the evolution of traffic state, but also account for the CF dynamics, while accurately capturing LC behaviors for multiple lanes.

## 5. Case Study

In this section, we first provide a brief overview of the dataset that utilized in the current study, and introduce the indicators for the analysis of trajectory accuracy. Then we evaluate the effectiveness of the proposed method under different penetration rates of PVs, and compare it with two conventional trajectory reconstruction methods. Lastly, we summarize the estimation results of LC position and analyze potential factors that could impact the estimation accuracy.

### 5.1 Data Description and Performance Indicators

The proposed method was evaluated using Next Generation Simulation (NGSIM, 2006) data, which is a publicly available dataset of human-driven vehicle trajectories and has been widely utilized for validating trajectories (Li et al., 2020). The dataset includes comprehensive trajectories of all vehicles recorded on both highways and arterials. The analysis focused on vehicle trajectory data collected from the leftmost lane (Lane a) and its adjacent lane (Lane b) of the southbound direction of U.S. Highway 101 in Los Angeles, California, spanning from 7:50 a.m. to 8:05 a.m. Fig. 8(a) depicts the study area and the upstream and downstream locations of virtual fixed sensors. The grey lines in Fig. 8(b) denote the ground-truth trajectories, while the blue dots indicate the detected points obtained from fixed sensors.

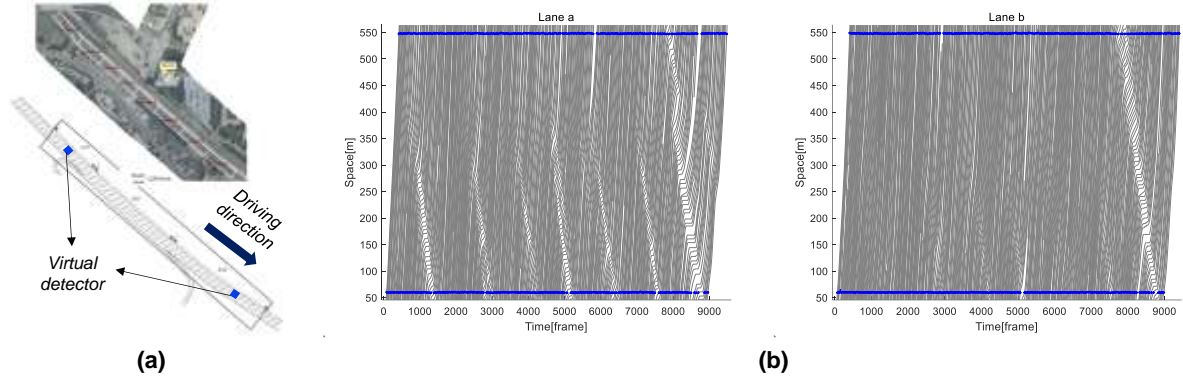

**Fig. 8. Study area of U.S. Highway 101: (a) Locations of virtual fixed sensors; (b) Ground-truth trajectories and detected points of two neighboring lanes.**

As the reconstructed outcomes were individual vehicle trajectories, certain macroscopic indicators like travel time and delay were deemed unsuitable. Therefore, three microscopic indicators, namely the Mean Absolute Error (MAE), Mean Absolute Percentage Error (MAPE), and Root Mean Square Error (RMSE), were employed to provide a comprehensive assessment of the proposed approach. These three indicators have been extensively utilized to assess reconstruction precision (Chen et al., 2021; Wang et al., 2020), by comparing the ground-truth and reconstructed trajectories on a per-second basis.

$$MAE = \frac{1}{T}\sum_{t=1}^{T} |X_{re}^t - X_{gt}^t| \tag{19}$$

$$MAPE = \left(\frac{1}{T}\sum_{t=1}^{T} \frac{|X_{re}^t - X_{gt}^t|}{X_{gt}^t}\right) * 100 \tag{20}$$

$$RMSE = \sqrt{\frac{1}{T}\sum_{t=1}^{T}(X_{re}^t - X_{gt}^t)^2} \qquad (21)$$

where $X_{re}^t$ and $X_{gt}^t$ denote the positions of ground-truth and the reconstructed trajectories at time $t$, respectively. MAE illustrates the average location error, while RMSE is more sensitive to extreme errors, and MAPE expresses the percentage of error relative to the length of the reconstructed trajectory.

Table 1. Parameters of Algorithms

| Variable | Description | Value |
|---|---|---|
| $\sigma$ | Width of spatial smoothing | 6 $m$ |
| $\tau$ | Width of temporal smoothing | 2 $s$ |
| $v_{free}$ | Propagation velocity of free flow | 24 $m/s$ |
| $v_{cong}$ | Propagation velocity of congested flow | -5 $m/s$ |
| $\hat{v}$ | Speed threshold between flow and congested flow | 15 $m/s$ |
| $\Delta v$ | Transition width around $\hat{v}$ | 3.6 $m/s$ |

Moreover, prior to the application of the trajectory reconstruction approach, the fundamental traffic flow parameters must be determined. Given the extensive literature on the study area (U.S. Highway 101), the calibrated parameters from previous research (Tsanakas et al., 2022) were adopted directly, as presented in Table 1. In order to ensure the desired level of detail in the velocity contour map, the spatial and temporal smoothing widths were set to low values of 6 m and 2 s, respectively.

## 5.2 General Results

To better illustrate the impact of the integrated framework, the proposed method was compared against two representative non-integrated approaches. The first approach is the micro-based method, in which trajectories were generated based on trajectory estimation algorithms without the fusion process. The second approach was named macro-based method that involved the reconstruction of trajectories directly from the velocity contour map. Specifically, the macro velocities were treated as individual vehicle velocities, and virtual trajectories were then generated based on kinematic equations. Macro-based methods have been widely used in previous trajectory reconstruction studies, and readers can refer to Jiang et al. (2017) and Tsanakas et al. (2022) for more details.

Three different penetration rates of PVs: 5%, 10%, and 15% were tested. The reconstructed trajectories under 10% penetration rate were displayed in Fig. 9. The red lines and grey dash lines represented the trajectories of PVs' and non-probe vehicles, respectively. The reconstructed trajectories are represented by color gradient lines. In this representation, the color of the lines corresponds to the error values, where lighter colors indicate higher error magnitudes. Ground-truth LC positions are indicated by orange circles, while the estimated LC positions are represented by green circles. Fig. 9(a)~(c) illustrated the performance of three methods on two adjacent lanes.

As shown in Fig. 9(b) with magenta rectangle, the trajectories reconstructed using micro-based methods pose two issues: 1) error accumulation is considerable as the number of estimated trajectories increases; 2) the trajectory estimated based on the leading PV overlaps with the following PV. The macro-based method alleviated the above issues, but the reconstructed trajectories remained too crowded because of the absence of LC estimation. Besides, the spatial distribution of these trajectories was not precise enough since the macro-based model ignored micro CF behaviors, as outlined by magenta rectangle in Fig. 9(c). In contrast, the trajectories became more plausible, and the aforementioned conflicts were mitigated after applying trajectory fusion algorithm. Furthermore, the two-stage fusion algorithm effectively captured the LC behaviors, and the estimated LC positions were in close agreement with the ground-truth positions, as indicated by the orange and green circles in Fig. 9(a). The overall comparison of three methods under different penetration rates was summarized in Table 2.

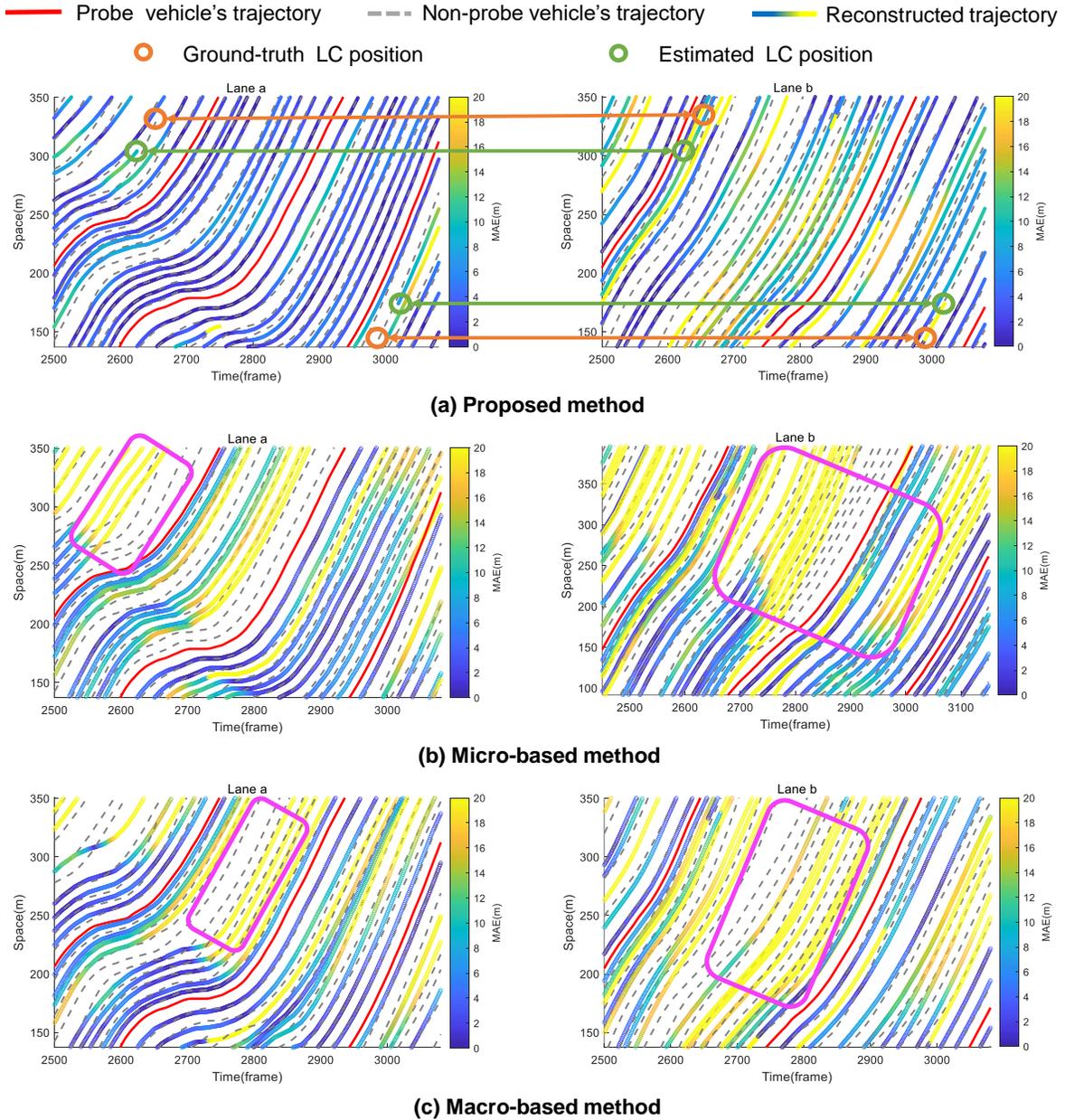

**Fig. 9. Reconstructed trajectories using three methods under 10% penetration rate of PVs: (a) Proposed method; (b) Micro-based method; (c) Macro-based method.**

During the 15-minute recording period, there were approximately 850 vehicles on the tested two adjacent lanes, and the reconstructed length was around 500 meters. Upon review of Table 2, it is evident that the penetration rate of PVs has a significant impact on both the non-integrated and proposed methods. In particular, for the proposed method, the values of MAE, MAPE, and RMSE decreased 49.80%, 50.87%, and 49.12%, respectively, as the penetration rate increased from 5% to 15%. This outcome can be attributed to the fact that a higher penetration rate of PVs can provide more accurate inputs and observed data to the integrated framework. Similarly, these indicators also dropped by around 50% for the other two compared methods.

Table 2. Comparison of three methods using NGSIM data

| Penetration Rate of PVs | Reconstructed Method | MAE (m) | MAPE (%) | RMSE (m) |
|---|---|---|---|---|
| 5% | **Proposed method** | **7.57** | **1.73** | **9.04** |
|  | Micro-based method | 17.27 (+128.14%) | 3.77 (+117.92%) | 23.62 (+161.28%) |
|  | Macro-based method | 11.44 (+51.12%) | 2.58 (+49.13%) | 14.56 (+61.06%) |
| 10% | **Proposed method** | **4.90** | **1.08** | **5.95** |
|  | Micro-based method | 11.91 (+143.06%) | 2.63 (+143.52%) | 16.05 (+169.75%) |
|  | Macro-based method | 8.34 (+70.20%) | 1.92 (+77.78%) | 10.55 (+77.31%) |
| 15% | **Proposed method** | **3.80** | **0.85** | **4.60** |
|  | Micro-based method | 8.10 (+113.16%) | 1.82 (+114.12%) | 10.50 (+128.26%) |
|  | Macro-based method | 6.56 (+72.63%) | 1.53 (+80.00) | 8.31 (+80.65%) |

The results demonstrated that the proposed method outperformed the other two methods across various penetration rates. The numbers in parentheses represented the percentage increase of errors compared to the proposed method. As presented in Table 2, MAE, MAPE and RMSE of the micro-based method under 5% penetration rate increased by 128.14%, 117.92% and 161.28%, respectively. Likewise, the macro-based method yielded higher values for these three performance indicators by 51.12%, 49.13%, and 61.06%, respectively. Owing to the disregard of the influence of traffic state and LC behaviors, the micro-based method exhibited the least effectiveness under all penetration rates. In contrast, the macro-based method performed better due to the imposition of traffic state constraints. Nonetheless, the errors remained significantly larger, exceeding 50% when compared with the proposed method. Upon the integration of the macro traffic state and micro CF and LC behaviors, the reconstructed trajectories of the proposed method demonstrated a remarkable enhancement in accuracy and smoothness.

On the other hand, it can be observed that the RMSE was generally larger than the MAE, indicating that the errors were unevenly distributed in the reconstructed region. Specifically, the estimated trajectories owned higher accuracy near the entrance and exit, whereas the errors tended to increase in the middle of the road. This phenomenon could be attributed to the deployment of fixed sensors in the entrance and exit, which provided approximate positions of each non-probe vehicle. Furthermore, the MAPE was found to be lower than 2% for all penetration rates, indicating that the location errors were within an acceptable range when compared with the reconstruction length.

### 5.3 Analysis of LC Position Estimation

As analyzed above, accurate estimation of LC behaviors holds the potential to enhance safety, efficiency, and optimize traffic flow on freeways. In order to further verify the efficacy of the proposed method, an analysis of the discrepancies between the ground-truth and estimated LC positions was conducted in this subsection. We classified the precision of LC position estimation into three categories: well-matching, moderate-matching, and failed-matching. Fig. 10(a) depicted the well-matching scenario, wherein the distance between the estimated and ground-truth LC positions is less than 30 meters. Given the sparsity of data input and the 500-meter length of reconstruction, an error within 30 meters can be regarded as highly accurate. On the other hand, Fig. 10(b) illustrated the moderate-matching scenario, wherein the consistency of LC behaviors was ensured during the reconstruction process, albeit with an error exceeding 30 meters between the ground-truth and estimated LC positions. However, it is worth noting that the consistency of LC behaviors cannot always be guaranteed due to the variations of traffic states and the constraints of safety distance. Such variations leaded to discrepancies in the estimated LC positions between the adjacent lanes, as demonstrated in Fig. 10(c), which we referred to as the failed-matching scenario.

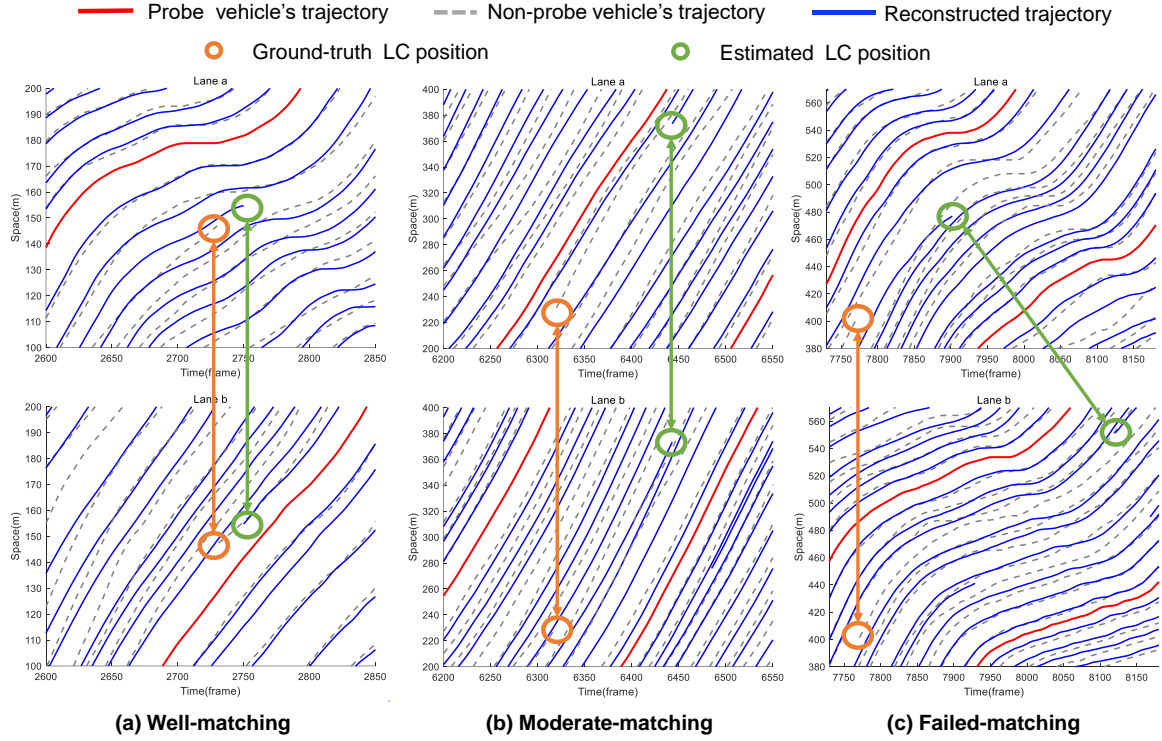

**Fig. 10.** Classification of the precision of LC position estimation: (a) Well-matching; (b) Moderate-matching; (c) Failed-matching.

Table 3 provided a summary of the three LC estimation scenarios under different penetration rates. The total number of ground-truth matched LC behaviors between Lane a and Lane b was recorded at 21 times, and the successful rate was calculated as the quotient of success matching scenarios and the total number of LC behaviors. As shown in Table 3, the proposed method exhibited effective LC position estimation for two-lane trajectory reconstruction, with a success rate exceeding 70% even under a 5% penetration rate. But owing to the highly sparse data environment, only a few instances of well-matching scenarios were observed. With a gradual increase in the penetration rate from 5% to 10%, the successful rate improved to 85% and well-matching scenarios doubled, which indicated that over 85% of the LC behaviors between adjacent lanes can be paired and approximately 30% of them are well-matched during the trajectory reconstruction process. Besides, no significant improvement was observed between the 10% and 15% penetration rates. Considering that LC positions were primarily determined by traffic state variations, the results suggested that the modified ASM had a good robustness in traffic state estimation even under low penetration rates.

**Table 3. Results of LC position estimation under different penetration rates.**

| Penetration Rate of PVs | Number of LC estimation scenarios | | | Successful rate |
|---|---|---|---|---|
| | Well-matching | Moderate-matching | Failed-matching | |
| 5% | 3 | 12 | 6 | 71.43% |
| 10% | 6 | 12 | 3 | 85.71% |
| 15% | 6 | 13 | 2 | 90.48% |

As seen from Fig. 10(c), the trajectory fusion algorithm at Stage 2 encountered failure in areas where significant changes in traffic state occurred. The main reason was that the driving heterogeneity made it more challenging to reconstruct high-accurate trajectories in such areas. Moreover, the error was cumulative and transitive, which further impacted the trajectory of subsequent vehicles. As a result, the estimated LC position may fail to satisfy the safety distance requirement in both lanes, leading to failed-matching scenarios. Therefore, considering the driving heterogeneity might improve the reconstruction accuracy and enhance the precision of LC position estimation, which will be one of the future research directions.

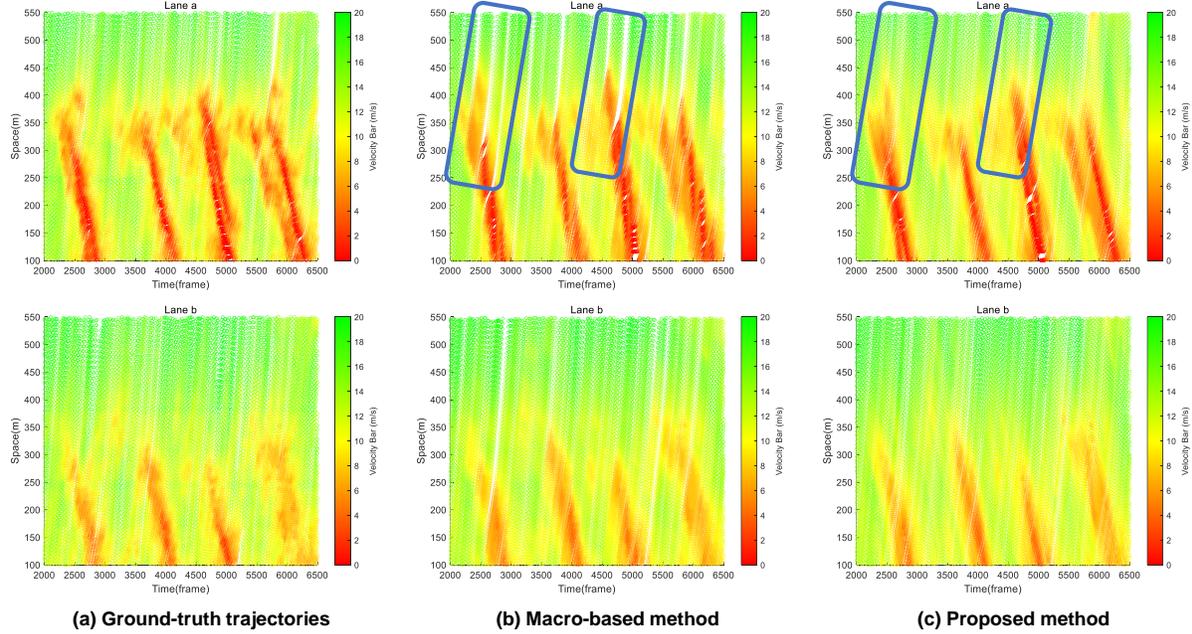

**Fig. 11. Comparisons of velocity contour map under 10% penetration rate: (a) Ground-truth trajectories; (b) Macro-based method; (c) Proposed method.**

In addition, to enable a comprehensive comparison between the macro-based and proposed methods, the velocity contour maps of ground-truth and reconstructed trajectories were generated for the entire region under the 10% penetration rate. These maps were presented in Fig. 11, wherein congested areas were depicted in red and free-flow areas in green. Notably, the proposed method demonstrated better consistency with the ground-truth trajectories, exhibiting minimal color discrepancies when compared to the macro-based method in both congested and free-flow areas, as indicated by dark blue rectangles. This observation underscores the accuracy and smoothness of the proposed method in reconstructing trajectories under diverse traffic conditions. Furthermore, the reconstructed trajectories successfully captured several congested traffic flows, which suggests the potential of the proposed method in enabling the analysis of more intricate microscopic traffic phenomena, such as emission estimation and traffic wave oscillations.

## 6. Conclusion and Future Work

This study presents a macro-micro framework for reconstructing fully sampled vehicle trajectories on multi-lane freeways using fixed and mobile sensors. The primary conclusions are highlighted below:

1) The proposed framework enables the trajectory reconstruction for multi-lane scenarios under the integration of macro and micro traffic models, enhancing the driving freedom of reconstructed trajectories and reproducing lane change behaviors.
2) Macroscopic velocity contour maps are established through fixed and mobile sensor data, which offer valuable reference for guiding LC positions and impose constraints on the collective movements of vehicles. Meanwhile microscopic trajectory estimation provides lane-based candidate trajectories for each non-probe vehicle.

3) The restriction presented by lane boundaries is effectively addressed through the utilization of the trajectory fusion algorithm, which solves the weights of candidate trajectories and identifies optimal LC positions across multiple lanes to jointly infer CF and LC behaviors.
4) Results show that the proposed method improves reconstruction accuracy, as evidenced by a substantial reduction in performance indicators relative to two non-integrated methods. Additionally, the success rate of LC position matching is over 85% even under a 10% penetration rate.

Nevertheless, this paper leaves several extensions for future research. Firstly, the trajectory fusion algorithm still faces challenges in matching lane change positions in several scenarios with rapid changes of traffic state. Hence, future research should ameliorate such failure scenarios, for example, by incorporating driver heterogeneity. This would further bolster the applicability of the proposed method in real-world traffic systems. Secondly, a small fraction of vehicles might execute multiple times of lane changes within the reconstruction range. However, the proposed method exhibits shortcomings in dealing with such vehicles as they are still detected by upstream and downstream sensors in the same lane. Exploring trajectory matching algorithms, such as filtering techniques (Wei et al., 2020), may present a promising solution to this challenge and constitute a valuable area for future research.

## Acknowledgement

This research is jointly sponsored by the National Natural Science Foundation of China (52125208), the Science and Technology Commission of Shanghai Municipality (No. 22dz1203200), and the Scientific and Technological Innovation 2030 - "New Generation Artificial Intelligence" Major Project (2022ZD0115602). The first author (X. Chen) is grateful to the China Scholarship Council (CSC) for financially supporting his visiting program at Tokyo Institute of Technology (No. 202206260138). The third author (T. Seo) is partially funded by a JSPS KAKENHI Grant-in-Aid for Scientific Research 20H02267.